\documentclass[11pt]{article}
\usepackage{graphicx}
\usepackage{amssymb}
\usepackage{amsmath}
\usepackage{amsthm}
\usepackage{amsfonts}
\usepackage{amssymb}
\usepackage{amssymb,color}

\def\'#1{\ifx#1i{\accent"13 \i}\else{\accent"13 #1}\fi}

\usepackage{url}

\newtheorem{theorem}{Theorem}
\newtheorem{lemma}{Lemma}

\newtheorem{remark}{Remark}

\newtheorem{problem}{Problem}

\title{On $ \omega \psi $-Perfect Graphs}

\author{G. Araujo-Pardo\footnotemark[2] \\
{\small {\tt garaujo@math.unam.mx}}
\and C. Rubio-Montiel\footnotemark[2] \footnotemark[3] \\
{\small {\tt christian.rubio@fmph.uniba.sk}}
}

\date{}

\begin{document}
\maketitle

\def\thefootnote{\fnsymbol{footnote}}
\footnotetext[2]{Instituto de Matem{\' a}ticas, Universidad Nacional Aut{\' o}noma de M{\' e}xico, 04510 M{\' e}xico City, Mexico.}
\footnotetext[3]{Department of Algebra, Comenius University, 842 48 Bratislava, Slovakia.}

\begin{abstract} 
In this paper, we generalize the concept of {\it{perfect graphs}} to other parameters related to graph vertex coloring. This idea was introduced by Christen and Selkow in 1979 and Yegnanarayanan in 2001. \\
Let $ a,b \in \{ \omega, \chi, \Gamma, \alpha, \psi \} $ where $ \omega $ is the clique number, $ \chi $ is the chromatic number, $ \Gamma $ is the Grundy number, $ \alpha $ is the achromatic number and $ \psi $ is the pseudoachromatic number. A graph $ G $ is \emph{$ ab $-perfect}, if for every induced subgraph $ H $ of $G$, $ a(H)$ equals $b(H) $. In this paper, we characterize the $ab$-perfect graphs when $a=\omega$ and $b=\psi$.
\end{abstract}
\section{Introduction}
Let $G$ be a simple graph. A vertex coloring $ \varsigma \colon V(G) \rightarrow \{1,\dots,k\} $ is called \emph{complete} if for each pair of distinct colors there exists an edge with these colors for its end vertices.
The \emph{pseudoachromatic number} $\psi(G)$ is the maximum $k$ for which a complete coloring of $G$ exists \cite{Gu}.

A vertex coloring of $ G $ is called \emph{proper} if  $\varsigma(u)\not=\varsigma(v)$ whenever $uv$ is an edge of $G$.
The \emph{achromatic number} $\alpha(G)$ is the maximum $k$ for which a proper and complete coloring of $G$ exists \cite{HHP}.

A vertex coloring of $ G $ is called \emph{Grundy}, if it is a proper vertex coloring having the property that for every two colors $ i $ and $ j $ with $ i<j $, every vertex colored $ j $ has a neighbor colored $ i $ (consequently, every Grundy coloring is a complete coloring).
The \emph{Grundy number} $\Gamma(G)$ is the maximum $k$ for which a Grundy coloring of $G$ exists \cite{Gr}. It is clear we have the following inequalities: 
\begin{equation}\label{des1}
\omega(G)\leq\chi(G)\leq\Gamma(G)\leq\alpha(G)\leq\psi(G)
\end{equation}
where $ \omega(G) $ and $\chi(G)$ denotes the clique number and the chromatic number of $G$ respectively.
\begin{figure}[!htbp]
\begin{center}
\includegraphics{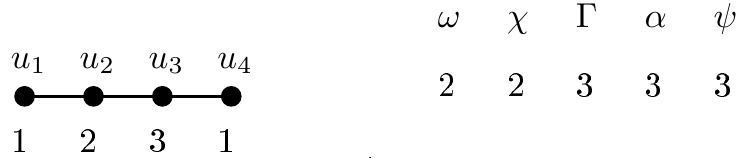}
\caption{\label{Fig3} A coloring of $P_{4}$ with $3$ colors.}
\end{center}
\end{figure}
A number of results of these invariants can be found in \cite{AMS,COTZ,CZ,Gu,HHP,YBS} and in the references therein. An interesting work related with the chromatic number, the achromatic number and the pseudoachromatic number was published by Yegnanarayanan in 2001 \cite{Y}. In that paper the connections between these three parameters with some combinatorial and computational topics are studied; such as the computational complexity, partitions, projective plane geometry and perfect graphs among others.

In particular, the authors of this paper and others have worked with the achromatic index and pseudoachromatic index of complete graphs vis a vis projective plane geometry (see \cite{AMS,AMRS}).

A graph $G$ without an induced subgraph $H$ of it is called $H$-\emph{free}. A graph $H_1$-free, $H_2$-free, $\dots$, $H_m$-free is called $(H_1,H_2,\dots,H_m)$-free.

Let $ a,b \in \{ \omega, \chi, \Gamma, \alpha, \psi \} $. A graph $ G $ is called \emph{$ ab $-perfect}, if for every induced subgraph $ H $ of $ G $ $ a(H)=b(H) $. This definition extends the usual notion of \textit{perfect graph} introduced by Berge \cite{B} in 1961; with this notation a perfect graph is denoted by $ \omega \chi $-perfect. We give some results about  this topic which we will use throughout this paper: L\'ovasz \cite{L} proved in 1972 that a graph is $ \omega \chi $-perfect, if and only if its complement is also $ \omega\chi $-perfect --the \emph{Weak Perfect Graph Theorem}--. Chudnovsky, Robertson, Seymour and Thomas \cite{CRST} proved in 2006 that $ G $ is $ \omega \chi $-perfect, if and only if $ G $ and $ G^c $ are $C_{2k+1}$-free for all $k\geq2$ --the \emph{Strong Perfect Graph Theorem}-- and Seinsche \cite{S} proved in 1974 that if a graph $ G $ is $ P_4 $-free then $ G $ is $ \omega \chi $-perfect.

The concept of the $ ab $-perfect graphs was introduced by Christen and Selkow in 1979 \cite{CS}. They considered $ ab $-perfect graphs for $ a,b\in \{ \omega, \chi, \Gamma, \alpha\}$, and they proved the following (see also \cite{GL}):
\begin{theorem}[Christen and Selkow \cite{CS}]\label{Grundy}
For any graph $ G $ the following are equivalent: {$\left\langle 1 \right\rangle $} $G$ is $\omega\Gamma$-perfect, {$\left\langle 2 \right\rangle $} $G$ is $\chi\Gamma$-perfect, and {$\left\langle 3 \right\rangle $}  $G$ is $P_4$-free.
\end{theorem}
\begin{theorem}[Christen and Selkow \cite{CS}]\label{achro}
For any graph $ G $ the following are equivalent: {$\left\langle 1 \right\rangle $} $G$ is $\omega\alpha$-perfect, {$\left\langle 2 \right\rangle $} $G$ is $\chi\alpha$-perfect, and {$\left\langle 3 \right\rangle $}  $G$ is $(P_4,P_3\cup K_2,3K_2)$-free.
\end{theorem}

Yegnanarayanan in \cite{Y}  also worked on this concept for $ a,b\in \{ \omega, \chi, \alpha, \psi \} $. In that work he stated:
\begin{enumerate}
\item For any finite graph $ G $ the following are equivalent: $\left\langle 1 \right\rangle $ $G$ is $\omega\psi$-perfect, $\left\langle 2 \right\rangle $ $G$ is $\chi\psi$-perfect, $\left\langle 3 \right\rangle $ $G$ is $\alpha\psi$-perfect and $\left\langle 4 \right\rangle $ $G$ is $C_{4}$-free.
\item Every $ \alpha \psi $-perfect graph is $ \chi \alpha $-perfect.
\item Every $ \chi \alpha $-perfect graph is $ \omega \chi $-perfect.
\end{enumerate}
Unfortunately Statement 1. is false (a counterexample is $P_4$ because it is $ C_4 $-free, but not $ \omega \psi $-perfect, see Figure \ref{Fig3}); i.e., $\left\langle 4 \right\rangle$ does not necessarily imply $\left\langle 1 \right\rangle$. Statement 2. is also false (then again, the counterexample is $P_4$, see Figure \ref{Fig3}) while Statement 3. is true (see \cite{CS}), but it is not well founded. Note that if a graph $ G $ is $ \omega \psi $-perfect then it is immediately $ \omega \chi $-perfect (see Equation \ref{des1}). That is, $ G $ is perfect in the usual sense and it is well known that $ G $ does not allow odd cycles (except triangles), or their complements (by the Strong Perfect Graph Theorem), therefore, the condition $ C_4 $-free is insufficient. With the aim of finding the correctness of this theorem we obtained the following result:
\begin{theorem}
\label{teo4}
For any graph $ G $ the following are equivalent:
\begin{description}
\item [{$\left\langle 1 \right\rangle $}] $G$ is $\omega\psi$-perfect,
\item [{$\left\langle 2 \right\rangle $}] $G$ is $\chi\psi$-perfect,
\item [{$\left\langle 3 \right\rangle $}] $G$ is $(C_4,P_4,P_3 \cup K_2,3K_2)$-free,
\item [{$\left\langle 4 \right\rangle $}] $G$ has one of the following structures:

i) If $G$ is a disconnected graph then: $G$ is an empty graph, $G$ has exactly two non-trivial complete graphs as components or $G$ has only one non-trivial $(C_4,P_4,P_3 \cup K_2,3K_2)$-free component.

ii) If $G$ is a connected graph then: $G$ is a complete graph or $G$ is the join between a complete graph and a disconnected graph $G'$ where $G'$ has exactly the same three possibilities as in the disconnected case. 
\end{description}
\end{theorem}
This paper is organized as follows: In the next section we study $ab$-perfect graphs for various graph parameters discussed above. In Section \ref{Results} we prove Theorem \ref{teo4}. Finally, in Section \ref{Open} we give some conclusions and open problems.
\section{Preliminaries}\label{pre}
Harary, Hedetniemi and Prins \cite{HHP} proved that for any graph $ G $ and for every integer $ a $ with $ \chi (G) \leq a \leq \alpha (G) $ there is a complete and proper coloring of $ G $ with $ a $ colors. Christen and Selkow \cite{CS} proved that for any graph $G$ and for every integer $b$ with $\chi(G)\leq b\leq \Gamma(G)$ there is a Grundy coloring of $G$ with $b$ colors (see also \cite{CZ} pg. 349).

Yegnanarayanan, Balakrishnan and Sampathkunar \cite{YBS} proved that, if $ 2\leq a \leq b \leq c $  there exists a graph $ G $ with chromatic number $ a $, achromatic number $ b $, and pseudoachromatic number $ c $. Chartrand, Okamoto, Tuza and Zhang \cite{COTZ} proved that for integers $a$, $b$ and $c$ with $2\leq a\leq b\leq c $ there exists a connected graph $G$ with $\chi(G)=a$, $\Gamma(G)=b$ and $\alpha(G)=c$, if and only if $a=b=c=2$ or $b\geq 3$ (see \cite{CZ} pg. 331 and 352).
\begin{figure}[!htbp]
\begin{center}
\includegraphics{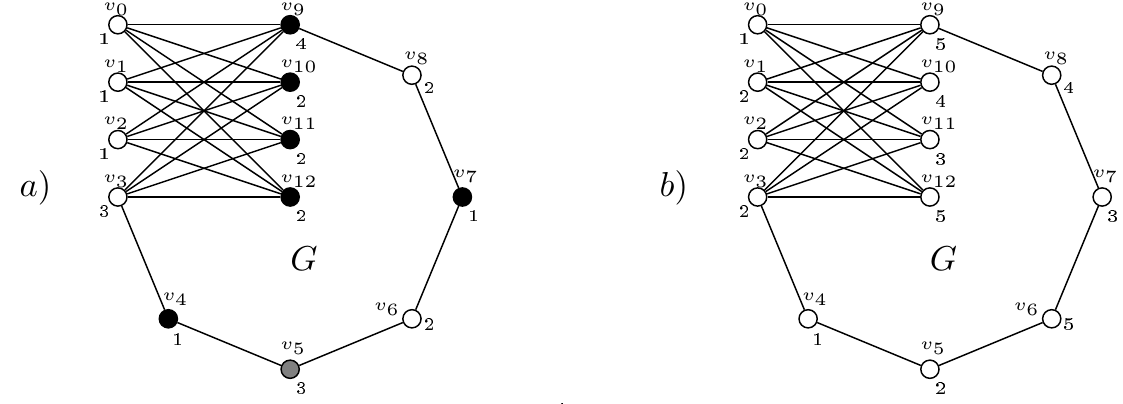}
\caption{\label{Fig1} Different coloration of $ G $.}
\end{center}
\end{figure}

Figure \ref{Fig1} shows a graph $ G $ where all of these parameters are different. In particular, $G$ has $\omega(G)=2$, $\chi(G)=3$, $ \Gamma(G)=4 $, $\alpha(G)=5$ and $\psi(G)=6$. Clearly, $ G $ does not contain $ C_3 $ as an induced graph because it is constructed identifying a $K_{4,4}$ with a $C_7$ cycle by an edge, as a result  $\omega(G)=2$. Figure \ref{Fig1}a) shows a proper vertex coloring of $ G $ with three colors (black, white and gray) then $ \chi(G)\leq3$ and $ G $ contains a $ C_7 $, therefore, $ \chi (G)=3 $. Also to the left, we show a Grundy vertex coloring of $ G $ with four colors (numbers) and $ \Gamma (G)\geq 4 $. Figure \ref{Fig1}b) shows a complete and proper vertex coloring of $ G $ with five colors and $ \alpha (G)\geq 5 $. It is not difficult to verify that $ \Gamma (G)\leq 4 $, $ \alpha (G) \leq 5 $ and $ \psi (G) = 6 $; details are omitted (see \cite{R} for details).

\begin{figure}[!htbp]
\begin{center}
\includegraphics{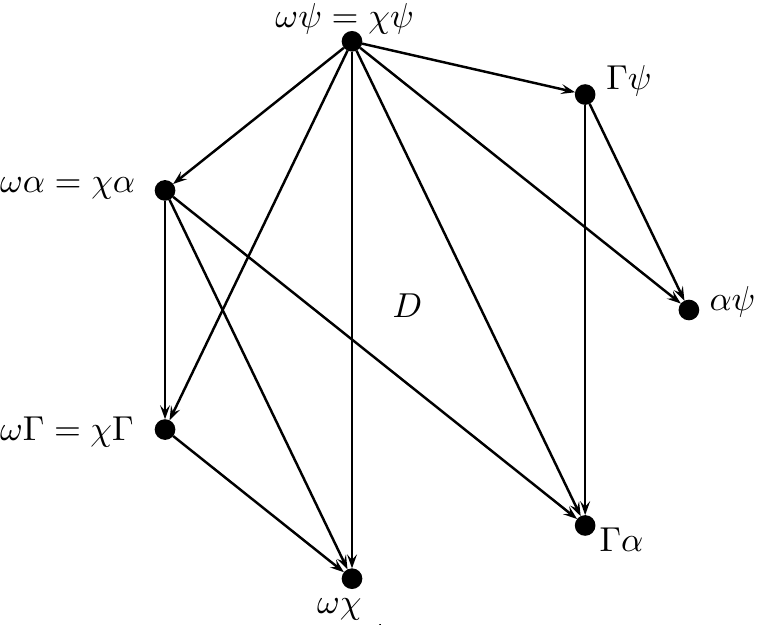}
\caption{\label{Fig2} Relationship between $ab$-perfect graphs with $a,b\in\{\omega,\chi,\Gamma,\alpha,\psi\}$.}
\end{center}
\end{figure}

In Figure \ref{Fig2} we show a directed transitive graph $ D $ where the vertices of $ D $ represent the classes of $ ab $-perfect graphs and their label is $ab$ respectively for $ a,b \in \{\omega,\chi,\Gamma,\alpha,\psi\} $. If two classes are equal they define the same vertex. If a class of $ab$-perfect graphs is contained in a class of $cd$-perfect graphs there exists an arrow from vertex $ab$ to vertex $cd$. To prove that $ D $ does not contain another arrow we have the following counterexamples: $ P_4 $, $ C_4 $, $ C_5 $ and $ P_3 \cup K_2 $ (see \cite{AR13,R} for details).

\section{On Theorem  \ref{teo4}} \label{Results}
To begin with, we establish some notation. For $A \subseteq V(G), $ let $\left\langle A\right\rangle $ be the subgraph of $G$ induced by $A$. Moreover, a graph $H$ is an induced subgraph of a $G$, briefly denoted by $H \leq G$, if there is a set $A \subseteq V(G)$, such that $\left\langle A\right\rangle $ is isomorphic to $H$. If we consider a set $ \{ v_1,\dots,v_n\} \subseteq V $ then we define $ \left\langle v_1,\dots, v_n \right\rangle := \left\langle \{v_1,\dots,v_n\} \right\rangle $. The set of vertices adjacent to the vertex $x$ is denoted by $N_G(x)$ (or only $N(x)$) and $N_G[x]=N(x)\cup\{x\}$ (or only $N[x]$). The maximum degree of vertices of $ G $ is denoted by $ \Delta(G) $ and $deg_G(v) $, or simply $ deg(v) $, denotes the degree of a vertex $ v $ in $ G $.
It is not difficult to note the following remarks:
\begin{remark} \label{R1} If $ G $ is a connected $ P_4 $-free graph then $ diam(G)\leq2  $.
\end{remark}
By the hereditary property of the $H$-free graphs we have the following:
\begin{remark} \label{R2} If $ G $ is an $ H $-free graph then every induced subgraph of $ G $ is also $ H $-free.
\end{remark}
The following lemma was proved in \cite{W2}, which follows the well-known Theorem \ref{teo3}.
\begin{lemma} \label{lemma1} If $ G $ is a connected $(C_4,P_4)$-free graph of order $ n $ then $ \Delta(G)=n-1 $.
\end{lemma}
\begin{theorem} \label{teo3}
A connected graph $ G $ is $ (C_4,P_4) $-free, if and only if $G$ is either a complete graph or there exists a set of connected $(C_4,P_4)$-free graphs $ \{ G_1,\dots,G_k \} $ for some $ k\geq2$ and $ m\geq1$, such that $G=K_{m}\oplus\overset{k}{\underset{i=1}{\cup}}G_{i}$.
\end{theorem}
The graphs $(C_4,P_4)$-free are also called \emph{trivially perfect graphs} (see \cite{Go1,Go2}), \emph{comparability graphs of trees} (see \cite{W1,W2}), or \emph{quasi-threshold graphs} (see \cite{MWW,YCC}).

Now we prove two short and useful lemmas:
\begin{lemma} \label{lemma2} If $H$ is one of the following three graphs: An empty graph, a graph with at most one non-trivial complete component, or a graph with at most two of them then  $\omega(H)=\psi(H)$.

\begin{proof} If $H$ is an empty graph or if $H$ has exactly one non trivial complete component then clearly $\omega(H)=\psi(H)$. We suppose that $H$ has two non-trivial complete components $K_{m_1}$ and $K_{m_2}$. Then $ \omega(H)=\max\{m_1,m_2\} $. Now we suppose $ \max\{m_1,m_2\}< \psi(H) $ and let $ \varsigma \colon V(H)\rightarrow\{1,\dots,\psi(H)\} $ a complete coloring of $H$. Note that each color class has at least two vertices in $K_{m_1}\cup K_{m_2}$, otherwise, as $ deg(v)\leq \max\{m_1,m_2\}-1 $ for any $ v \in V(H) $ then $\varsigma$ will not be complete, consequently $ \psi(H)\leq \frac{m_1+m_2}{2} \leq \max\{m_1,m_2\} $, a contradiction and hence $ \psi (H)=\max\{m_1,m_2\} $.
\end{proof}
\end{lemma}

\noindent{\bf Proof of Theorem \ref{teo4}.}

The proof of $\left\langle 1 \right\rangle \Rightarrow \left\langle 2 \right\rangle$ immediately follows from (\ref{des1}).

To prove $\left\langle 2 \right\rangle \Rightarrow \left\langle 3 \right\rangle$ note that, if $ H \in \{ C_4, P_4, P_3 \cup K_2, 3K_2 \} $ then $ \chi (H)=2 $, and $ \psi (H)=3 $, hence the implication is true.

We consider two cases to prove $\left\langle 3 \right\rangle \Rightarrow \left\langle 4 \right\rangle$. First, if $G$ is a disconnected graph then $G$ should be one of: $G$ is an empty graph, $G$ has exactly two non-trivial complete graphs as components or $G$ has only one non-trivial $(C_4,P_4,P_3 \cup K_2,3K_2)$-free component since $G$ is $(P_3\cup K_2,3K_2)$-free. From this point onward, if $G$ is one of these three disconnected graphs we will say that $G$ has a \emph{special structure}. \\

Second, if $G$ is a connected graph, since $G$ is $(C_4,P_4)$-free by Theorem \ref{teo3}, $G$ is either a complete graph or there exists a set of connected $(C_4,P_4)$-free graphs $ \{ G_1,\dots,G_k \} $ for some $ k\geq2$ and $ m\geq1$ such that $G=K_{m}\oplus\overset{k}{\underset{i=1}{\cup}}G_{i}$. Let $G'=G\setminus K_m$, then, $G'$ has at most two non-trivial components since it is $3K_2$-free, as $G'$ is also $P_3\cup K_2$-free we conclude $G'$ has a special structure.

To prove $\left\langle 4 \right\rangle \Rightarrow \left\langle 1 \right\rangle$ suppose that $G$ is a disconnected graph, such that $G$ is either an empty graph or has two non-trivial complete components. As a consequence of  Remark \ref{R2} and Lemma \ref{lemma2}, $G$ is $\omega\psi$-perfect. If $G$ has exactly one non-trivial $(C_4 , P_4 , P_3\cup K_2 , 3K_2 )$-free component $G'$ we conclude, by induction over $ \left| V(G) \right| $, that $G'$ is $\omega\psi$-perfect, and also $G$ since $\omega(G)=\omega(G')$ and $\psi(G)=\psi(G')$.  \\

Now we suppose that $G$ is a connected graph. If $G$ is a complete graph then it is $\omega\psi$-perfect because every induced subgraph of $G$ is a complete graph (it is also a consequence of Remark \ref{R2} and Lemma \ref{lemma2}). Finally, if $G$ is the join between a complete graph and a disconected graph $G'$ with the special structure we will divide the proof into two cases and in both we use induction over $ \left| V(G) \right| $: \\

{\bf Case 1)} If $G'$ is one of: An empty graph, a graph with at most one non-trivial complete component, or a graph with at most two of them, we can write $G=K_{n_{1}}\oplus\left(K_{n_{2}}\cup K_{n_{3}}\cup n_{4}K_{1}\right) $ where $n_1,n_2,n_3\geq1$ and $n_4\geq0$, if $n_4=0$ it means $n_4K_1=\emptyset$.\\
If $ \left| V(G) \right|=3 $ then $ G \cong P_3 $ is $ \omega \psi $-perfect. Assume that, if $4\leq  \left| V(G) \right| \leq m $ then $ G $ is $ \omega\psi $-perfect, and let $ G=K_{n_{1}}\oplus\left(K_{n_{2}}\cup K_{n_{3}}\cup n_{4}K_{1}\right) $ be a graph, such that $ n_1+n_2+n_3+n_4=m+1 $. If $ H $ is a connected induced subgraph of $ G $ with $ \left| V(H) \right| < \left| V(G) \right| $ then $ H=K_{m_{1}}\oplus\left(K_{m_{2}}\cup K_{m_{3}}\cup m_{4}K_{1}\right) $ where $ m_i \leq n_i $ for all $ i \in \{ 1,\dots,4 \} $ and $ m_1+m_2+m_3+m_4 \leq m $, by induction hypothesis $ H $ is $ \omega\psi $-perfect, that is $ \omega (H)= \psi (H) $. \\
If $ H $ is a disconnected induced subgraph of $ G $, such that $ \left| V(H) \right| < \left| V(G) \right| $ then $ H=K_{m_1}\cup K_{m_2}\cup m_3K_1  $ where $ m_1 \in \mathbb{Z}^+ $ and $ m_2, m_3 \in \mathbb{N} $, applying Lemma \ref{lemma2}  we conclude $\omega (H)=\psi(H) $. \\
Futhermore, we will prove $ \omega (G) = \psi (G) $: Let $ n=\max\{n_{2},n_{3}\}, $ clearly $ \omega(G)=n_1+n $. Suppose that $ \psi(G)>n_{1}+n $, and let $ \varsigma \colon V(G)\rightarrow\{1,\dots,\psi(G)\} $ be a complete coloring of $ G $, first of all we also suppose that $ n_4>0 $ and consider $ u $ a vertex of $ n_4K_1 $; as the degree of $ u $ is $ n_1 $ $ \varsigma(u) $ meets at most $ n_1 $ different chromatic classes from $ u $. In consequence there must exists another vertex $ v $, such that $ \varsigma(u) = \varsigma (v) $. Clearly, $ N(u) \subseteq N[v] $ then $ \psi(G)=\psi(G\setminus u) $, but $ \omega (G) = \omega (G\setminus u) $ and $\omega (G\setminus u) < \psi(G\setminus u)$ which is not possible by the induction hypothesis ( $G\setminus u $ is an induced subgraph of $ G $ and it is $\omega\psi$-perfect). Then, if $ \psi (G) > n_1 +n $, $ G=K_{n_{1}}\oplus\left(K_{n_{2}}\cup K_{n_{3}}\right) $ and $ n = \max \{ n_2,n_3 \} $. Let $ u \in V(K_{n_2}\cup K_{n_3}) $, such that $ \varsigma(u) \notin \{\varsigma(v) : v \in V(K_{n_1}) \} $. As the neighbors of $ u $ are at most $ n_1 +n -1 $ the chromatic class of $u$ meets, from $u$, at most $ n_1+n-1 $ different chromatic classes and there necessarily exists another vertex $ v $ in $ V(K_{n_2}\cup K_{n_3}) $, such that $ \varsigma(u) = \varsigma (v) $ and $ \psi(G) \leq n_{1}+\frac{n_2+n_{3}}{2} \leq n_1 + n $, which is a contradiction. Consequently, $ \omega(G)=\psi(G) $, hence, $ G $ is $ \omega\psi $-perfect.\\
 
{\bf Case 2)} If $G'$ has exactly one non-complete $(C_4 , P_4 , P_3\cup K_2 , 3K_2 )$-free component $H'$ we can write $G=K_{n_{1}}\oplus\left(H'\cup n_{2}K_{1}\right)$ where $n_1,n_2\geq1$.\\
If $ \left| V(G) \right|=5 $ then $ G=K_1 \oplus (P_3\cup K_1)$. It is not difficult to see that it is $ \omega \psi $-perfect. Assume that any graph $ G $ with this structure and order $6\leq  \left| V(G) \right| \leq m $ is $ \omega\psi $-perfect, and let $ G=K_{n_{1}}\oplus\left(H'\cup n_{2}K_{1}\right) $, such that $ n_1+n'+n_2=m+1 $ where $ n'= \left| V(H') \right|$. If $ H $ is a connected induced subgraph of $ G $, such that $ \left| V(H) \right| < \left| V(G) \right| $, then by Remark \ref{R2} $ H $ is $(C_4,P_4,P_3\cup K_2,3K_2)$-free. Therefore, by $\left\langle 3 \right\rangle \Rightarrow \left\langle 4 \right\rangle$ we have that $H$ is a complete graph or $H$ is the join between a complete graph and a disconected graph $H''$ with the special structure (*). In any case $H$ is $ \omega \psi $-perfect; using the previous case or the induction hypothesis, in particular $ \omega(H) = \psi (H) $ (**). Furthermore, if $ H $ is an induced disconnected subgraph of $ G=K_{n_1}\oplus(H'\cup n_2K_1) $, such that $ \left| V(H) \right| < \left| V(G) \right| $ then $ H $ is an induced disconnected subgraph of $ H'\cup n_{2}K_{1} $, if $ H \leq n_2K_1 $ we apply Lemma \ref{lemma2} and obtain $ \omega (H) = \psi (H) $.  On the other hand, if $ H \not \leq n_2K_1 $, $ \psi (H) = \psi (H\setminus n_2K_1) $ and $ \omega (H) = \omega (H\setminus n_2K_1) $. As $ H\setminus n_2K_1 \leq H' $ and $ H' $ is a proper induced connected subgraph of $G$, for (*) and (**) $ H' $ is $ \omega\psi $-perfect, in particular $ \omega (H\setminus n_2K_1)=\psi (H\setminus n_2K_1) $ and $ \omega (H)=\psi (H) $.\\
Finally, to prove that $ \omega (G) = \psi (G) $ we note that $ n_2>0 $ and using exactly the same argument as in the previous case (for $n_4>0$) we arrive to a contradiction. Therefore, $G$ is $ \omega\psi $-perfect.
{\ \rule{0.5em}{0.5em}\vskip 12pt}
\section{Open problems}\label{Open}
To conclude, we wish to present some problems that, in our opinion, might be of interest. The first one is a generalization of the problem only for the chromatic, achromatic and pseudoachromatic parameters as we state in Section \ref{pre} (see \cite{COTZ,CZ,YBS}), and the second one can be seen as a problem in graph theory using the language of perfectness.
\begin{problem}
Let $ 2 \leq a \leq b \leq c \leq d \leq e  $. Does there exist a graph $ G $, such that $ \omega(G)=a $, $ \chi(G)=b $, $ \Gamma(G)=c $, $ \alpha(G)=d $, and $ \psi (G)=e $?
\end{problem}

\begin{problem} \label{pr2}
Characterize the $\alpha\psi$-perfect graphs.
\end{problem}

Suppose we have $k$ graph invariants $\{\beta_1,\beta_2,\ldots,\beta_k\}$ such that every graph $G$ satisfies $\beta_1(G)\leq \beta_2(G) \leq \ldots \leq \beta_k(G)$ for any two different $i,j\in \{1,\ldots,k\}$. Let ${\cal{B}}_{ij}$ be the class of $\beta_i\beta_j$-perfect graphs. By definition of ${\cal{B}}_{ij}$ as a hereditary class, there exists a family ${\cal{F}}_{ij}$ of graphs that should not be contained in ${\cal{B}}_{ij}$, such that ${\cal{B}}_{ij}$ is exactly the class of ${\cal{F}}_{ij}$-free graphs. It is also clear that: 
$${\cal{B}}_{ij}={\cal{B}}_{i,i+1}\cap{\cal{B}}_{i+1,i+2},\ldots,\cap{\cal{B}}_{j-1,j}$$
Therefore, we are interested in determining ${\cal{F}}_{i,i+1}$ for each $i=\{1,\dots, k-1\}$.

We have to analyze the inequalities that appear in Equation \ref{des1}, as we noted in the introduction, the $\omega\chi$-perfect graphs are $C_{2k+1}$-free and  $C^c_{2k+1}$-free (see \cite{CRST}) and the $\chi\Gamma$-perfect graphs are $P_4$-free (see \cite{CS}). To characterize the $\Gamma\alpha$-perfect graphs a nonredundant family (in the sense that none is an induced subgraph of any other) of ten forbidden subgraphs for $\Gamma\alpha$-perfect graphs is exhibited in \cite{CS} which includes the graphs $P_3\cup K_2$ and $3K_2$. 

For the case of $\alpha\psi$-perfect graphs: First, Hedetniemi conjectured that for any tree $T$, $\alpha(T)=\psi(T)$, however, Edwards showed that the conjecture is false by exhibiting an infinite family of trees $T_i$ for which $\alpha(T_i)<\psi(T_i)$ (see \cite{E}). The smallest example $T_1$ has 408 vertices, nevertheless, it is not known if it is the smallest one. Therefore, there is at least one tree $T'$ in the family of forbidden graphs of $\alpha\psi$-perfect graphs. Second, it is known that $\alpha(C_n)=\psi(C_n)$ if and only if $n\not = 2x^{2}+x+1$ \cite{HM,Y2}, ie, a nonredundant infinite family of cycles is contained in the family of forbidden graphs for the class of $\alpha\psi$-perfect graphs which includes the graph $C_4$. 

\subsection*{\bf{Acknowledgment}}
The authors wish to thank the anonymous referees of this paper and also of \cite{AR13} for their kind help and valuable suggestions which led to an improvement of this paper. 

Research supported by CONACyT-M{\' e}xico under Projects 166306, 178395 and PAPIIT-M{\' e}xico under Project IN104915.

\end{document}